\documentclass[a4paper,12pt]{amsart}

\usepackage{epsfig}
\usepackage{graphicx}

\usepackage[latin1]{inputenc}
\usepackage[T1]{fontenc}
\usepackage{indentfirst}
\usepackage{amssymb}
\usepackage{eufrak}
\usepackage{amsmath}
\usepackage{amsfonts}
\usepackage{amsthm}
\usepackage{mathrsfs}
\usepackage[bookmarks]{hyperref}
\usepackage{xypic}
\newcommand{\R}{\mathop{\rm Re}\nolimits}

\newcommand{\Log}{\mathop{\rm Log}\nolimits}

\newcommand{\Arg}{\mathop{\rm Arg}\nolimits}
\newcommand{\Vol}{\mathop{\rm Vol}\nolimits}

\newcommand{\Area}{\mathop{\rm Area}\nolimits}
\newcommand{\Alga}{\mathop{\rm Alga}\nolimits}
\newcommand{\supp}{\mathop{\rm supp}\nolimits}
\newcommand{\Alg}{\mathop{\rm Alg}\nolimits}
\newcommand{\Critv}{\mathop{\rm Critv}\nolimits}

\input epsf

\def \square{\smallskip \hfill \vrule width 5 pt height 7
pt depth - 2 pt \smallskip }

\newenvironment{prooof}
{\noindent {{\it Proof} \;}}{\hspace*{\fill}\square\vskip 8pt}

\theoremstyle{plain}
\newtheorem{Lem}[subsection]{Lemma}
\newtheorem{The}[subsection]{Theorem}
\newtheorem{Cor}[subsection]{Corollary}
\newtheorem{Pro}[subsection]{Proposition}
\theoremstyle{definition}
\newtheorem{Rem}[subsection]{Remark}

\newtheorem{Def}[subsection]{Definition}

\begin{document}

\title[Coamoebas of Complex Algebraic Plane Curves and Gauss  Map]{Coamoebas of Complex Algebraic Plane Curves and The Logarithmic Gauss  Map}

\author{Mounir Nisse}
\address{Institut de Math{\'e}matiques de Jussieu (UMR 7586), Universit{\'e}
  Pierre et Marie Curie, Analyse Alg{\'e}brique\\ 175, rue du 
Chevaleret,\\ 75013
Paris} 
\email{nisse@math.jussieu.fr} 
\maketitle
%
%
%
%
%
\begin{abstract}
The coamoeba of any complex algebraic plane curve $V$ is its image in the real torus under the argument map.
The area counted with multiplicity of the coamoeba of any algebraic curve in $(\mathbb{C}^*)^2$ is bounded in terms of the degree of the curve. We show in this Note that up to multiplication by a constant in $(\mathbb{C}^*)^2$, the complex algebraic plane curves  whose coamoebas are of maximal area (counted with multiplicity) are defined over $\mathbb{R}$, and their real loci are Harnack curves possibly with ordinary real isolated double points  (c.f. \cite{MR-00}). In addition, we characterize the complex algebraic plane curves such that their coamoebas contain no extra-piece.
\end{abstract}

\setcounter{tocdepth}{1} \tableofcontents

\section{Introduction}

There is a very strong relationship between the amoeba and coamoeba of a complex algebraic hypersurface. Given an amoeba of a hypersurface, it is well known by Mikael Passare and Hans Rullg\aa rd \cite{PR-04} that the spine of the amoeba is equipped with a structure of tropical  hypersurface which is dual to some convex subdivision of the Newton polytope. We also know that amoeba reaches infinity by several tentacles \cite{V2-02}. The coamoeba covers the hyperplanes orthogonal to the external  edges of some subdivision of the Newton polytope by several {\em extra-pieces }. An extra-piece goes away from some  hyperplane orthogonal to some inner edge of the subdivision of the Newton polytope. These hyperplanes are equipped with a frames spanned by the normal vectors to the sides of the element of the subdivision containing those edges, and these frames determine the position of the coamoeba relatively to these hyperplanes.

The purpose of this Note is to describe the relations and the similarities which exist between amoebas and coamoebas of  complex algebraic plane curves.
For a polynomial in two variables $f\in \mathbb{C}[z,w]$, let $V$ be the curve of its zero set in $(\mathbb{C}^*)^2$ and $\Delta$ its Newton polygon, which is the convex hull in $\mathbb{R}^2$ of the index of the non zero coefficients of $f$. The {\em amoeba} $\mathscr{A} \subset \mathbb{R}^2$ of the curve $V$ is the image of $V$ by the logarithmic map $\Log$; where $\Log :(\mathbb{C}^*)^2\rightarrow\mathbb{R}^2$ is the map defined by : $\Log (z,w)=(\log \parallel z\parallel , \log \parallel w\parallel )$
( c.f. \cite{GKZ-94}).

The amoeba's complement has a finite number of convex connected components, corresponding to domains of convergence of the Laurent series expansions of the rational function $\frac{1}{f}$.
We know that the spine $\Gamma$ of the amoeba $\mathscr{A}$ has a structure of a tropical curve  in $\mathbb{R}^2$ (proved by M. Passare and H. Rullg\aa rd in 2000 \cite{PR-04}, and independently by G. Mikhalkin in 2000). In addition the spine of the amoeba is dual to some coherent (i.e. convex) subdivision $\tau$ of the integer convex polygon $\Delta$. It is shown by M. Forsberg, M. Passare and A Tsikh that the set of vertices of $\tau$ is in bijection with the set of complement components of $\mathscr{A}$  in $\mathbb{R}^2$ \cite{FPT-00}. It was shown by Mikael Passare and Hans Rullg\aa rd in \cite{PR-04} that the amoebas of complex algebraic plane curves have finite area, and the upper bound on the area is in terms of the the degree of the curve. In addition, it was shown by Grigory Mikhalkin and Hans Rullg\aa rd, that up to multiplication by a constant in $(\mathbb{C}^*)^2$ the complex 
algebraic plane curves whose amoebas are of maximal area  are defined over $\mathbb{R}$ and, furthermore their real locus are a Harnack curves possibly with ordinary real isolated double points  (c.f. \cite{MR-00}).
The {\em coamoeba}  denoted by $co\mathscr{A}\subset (S^1)^2$ of the curve $V$ is the image of $V$ under the argument mapping $\Arg : (z, w)\mapsto (e^{i\arg (z)}, e^{i\arg (w)})$. This  terminology has been introduced by M. Passare and A. Tsikh in 2000.
It is shown in \cite{N2-07} that the complement components of the closure in the flat torus of the coamoeba of a complex algebraic hypersurface of complex dimension $n-1$ and defined by a polynomial $f$ with Newton polytope $\Delta$ are convex and their number don't exceed $n!\Vol (\Delta )$.

In this  Note, we will prove a similar properties of  the complex algebraic plane curves coamoebas. More precisely we  
characterize the curves  such that the area of their coamoebs counted with multiplicity is maximal. Using the fact that the set of critical point $Crit_{\Arg}$ of the argument map is equal to that of the logarithmic map which is the inverse image by the logarithmic Gauss map of the real projective line, and some properties of extra-pieces introduced in \cite{N1-06} and \cite{N3-07}, we obtain the following results.

\vspace{0.3cm}

\begin{The} Let $V$ be a complex algebraic plane curve defined by a polynomial with Newton polygon $\Delta$. Then the following statements are equivalent:
\begin{itemize}
\item[(i)]  The cardinality of the   inverse image by the map $ pr: co\mathscr{A}\setminus \Arg (Crit_{\Arg})\rightarrow T\setminus pr\circ \Arg (Crit_{\Arg})$ of any point in the image is a constant integer $k$ which don't exceed  $2\Area (\Delta )$,
\item[(ii)] the coamoeba of $V$ has no extra-piece.
\end{itemize}
\end{The}

\newpage

\begin{The}
Let $V$ be a complex algebraic plane curve defined by a polynomial with Newton polygon $\Delta$. Then  the area counted with multiplicity of the coamoeba of $V$ cannot exceed $2 \pi^2\Area (\Delta )$, and we have 
the following equivalent statements:
\begin{itemize}
\item[(i)] $\Area_{mult}(co\mathscr{A}) = 2 \pi^2\Area (\Delta )$,
\item[(ii)] The curve $V$ is real up to multiplication by a constant in $\mathbb{C}^*$, and its real part $\mathbb{R}V$ is a Harnack curve possibly with ordinary real isolated double points.
\end{itemize}
\end{The}

\begin{Cor} Let us denote by $\Alga (V)$ the image in $T$ of $co\mathscr{A}\setminus \Arg (Crit_{\Arg})$ under the projection $pr$, and assume that the coamoeba of $V$ has no extra-piece. Then we have the equality $\Area_{mult}(co\mathscr{A}) = k \Area (\Alga (V))$. In particular, if $V$ is a Harnack curve, then $k=2\Area (\Delta )$ and $\Area (\Alga (V)) = \pi^2$.
\end{Cor}

\begin{Cor} Let $V$ be a complex algebraic plane curve defined by a polynomial with Newton polygon $\Delta$ such that $\Area (\Delta ) = \frac{1}{2}p$ where $p$ is a prime integer. Then the coamoeba of $V$ has no extra-piece if and only if $V$ is a Harnack curve or the spine of its amoeba has only one vertex and such that  the polynomial defining $V$ is maximally sparse.
\end{Cor}

\noindent The extra-pieces are subsets of the coamoeba with no-vanishing area and such that their boundary has a no discrete  component contained in the  set of critical values of the argument map. We can remark that if we assume that the closure in the torus of the coamoeba  is not all  the torus, then the no-linear components of the boundary of the coamoeba are contained in the boundary of some extra-pieces (for more details see \cite{N1-06} or \cite{N3-07}).
They play the role of the tentacles of the amoeba which correspond to the external edges of the subdivision $\tau$ dual to the spine of the amoeba. Except that extra-pieces correspond to some inner edges of the subdivision $\tau$.

\section{Preliminaries}

Let $V$ be an algebraic hypersurface in $(\mathbb{C}^*)^n$ defined by
the polynomial
$$
f(z) =\sum_{\alpha\in \supp (f)} a_{\alpha}z^{\alpha},\,\,\,\,\quad\quad\quad
z^{\alpha}=z_1^{\alpha_1}z_2^{\alpha_2}\ldots z_n^{\alpha_n},
\quad\quad\quad\quad\quad\quad\quad\quad  (1)
$$
where $a_{\alpha}$ are non-zero complex numbers and $\supp (f) = \{ \alpha\in\mathbb{Z}^n \mid a_{\alpha}\ne 0\}$ is
the support of $f$, and we
denote by $\Delta$ the Newton polytope of $f$ (i.e., the convex
hull in $\mathbb{R}^n$ of $\supp (f)$). We assume in this Note that the components of $\alpha$ are natural.

\subsection{Gauss map and the critical points of  the complex logarithmic map} ${}$
 
\vspace{0.2cm}

We use  here the definition of the Gauss map of an hypersurface in the
algebraic torus given by Kapranov in \cite{K1-91}.
\begin{Def}
Let $(\mathbb{C}^*)^n$ be the algebraic torus and $V\subset
(\mathbb{C}^*)^n$ be an algebraic hypersurface. For each
$z\in(\mathbb{C}^*)^n$ let $l_z :  h\longmapsto z.h$ be the translation
by $z$. The {\em Gauss map} of the hypersurface $V$ is the rational map
$\gamma : V \longrightarrow \mathbb{CP}^{n-1}$ taking a smooth point
$z\in V$ to the hyperplane $d(l_{z^{-1}})(T_zV)\subset
T_1((\mathbb{C}^*)^n)$, where $ \mathbb{CP}^{n-1} =
\mathbb{P}(T_1((\mathbb{C}^*)^n))$
\end{Def}

$$
\gamma  : V \longrightarrow  \mathbb{CP}^{n-1}
$$

$$
z \longmapsto T_zV\subset
T_z((\mathbb{C}^*)^n)\stackrel{t^*_{z^{-1}}}{\longrightarrow}
t^*_{z^{-1}}(T_zV)\subset T_1((\mathbb{C}^*)^n)
$$

$$
z\longmapsto [t^*_{z^{-1}}(T_zV)]
$$

\noindent where $[t^*_{z^{-1}}(T_zV)]$ is the class of $t^*_{z^{-1}}(T_zV)$ in $
\mathbb{CP}^{n-1}$.

We assume from  now that $V\subset (\mathbb{C}^*)^2$ is an algebraic curve
defined by a polynomial $f$ and nowhere singular. Let $z\in V$,
$U\subset (\mathbb{C}^*)^2$ be a small open neighborhood of $z$ and $\Log_U :
U\longrightarrow \mathbb{C}^2$ be a choice of a  branch of   the holomorphic logarithm function
 and then we take the image of $U\cap V$ by $\Log_U$. Hence the
vector $\gamma (z)\in \mathbb{C}^2$ ( looking at it as
$T_{\Log_U(z)}\mathbb{C}^2$) tangent to $\Log_U(U\cap V)$ at the point
  $\Log_U(z)$ is a complex vector, which is independent of the choice
  of the holomorphic branch. So it defines a map from $U\cap V$ (a
  rational map from  all $V$) to $\mathbb{CP}^1$ , called the
  {\em logarithmic Gauss map} of the curve $V$. Explicitly one have:

\begin{eqnarray*}
\gamma :&V&\longrightarrow\,\,\,\mathbb{CP}^1\\
&z&\longmapsto\,\,\, [\gamma (z)]
\end{eqnarray*}

\noindent where $[\gamma (z)]$ denotes the class of the vector $\gamma
(z)$ in $\mathbb{CP}^1$. We decompose the identity as follow:
$$
U\cap
V\stackrel{\Log_U}{\longrightarrow}\mathbb{C}^2\stackrel{\exp}{\longrightarrow}(\mathbb{C}^*)^2\stackrel{f}{\longrightarrow}\mathbb{C}
$$
Let us look now to the amoeba $\mathscr{A}$ of the curve $V$, and
more precisely to $\Log (U\cap V)$ (where here $\Log (z_1,z_2) =
(\log \parallel z_1\parallel , \log \parallel z_2\parallel )$). It's clear that $\Log
(U\cap V) = \R \circ \Log_U (U\cap V)$ where $\R (z_1,z_2) =
(\R (z_1),\R (z_2))$ and $\R$ denotes the real part. So if
$t=\Log_U(u)$ where $u = (u_1,u_2)\in U$, then $dt = \frac{du}{u}$, and we obtains
\begin{eqnarray*}
\gamma (z)& =& \frac{df(\exp (t))}{dt}_{\mid t=\Log_U(z)}\\
&=&(u_1\frac{\partial f}{\partial u_1}(u), u_2\frac{\partial f}{\partial u_2}(u))_{\mid u=z}\\
&=&(z_1\frac{\partial f}{\partial z_1}(z), z_2\frac{\partial f}{\partial z_2}(z))
\end{eqnarray*}
\noindent Hence a point $z\in V$ is a critical point of the function
$\Log = \R\circ\Log_U$ is equivalent to saying that the tangent space
of $\Log_U (U\cap V)$ at the point $\Log_U (z)$ contains a non zero
imaginary vector $v$, which means that $T_{\Log_U (z)}(\Log_U (U\cap
V))$ is the complexification of a real line and then it is invariant under complex conjugation so
 $\gamma (z) \in
\mathbb{RP}^1\subset \mathbb{CP}^1$.

\noindent We have the following Proposition and its proof given by Mikael Passare at IMJ, l'institut de Mathématiques de Jussieu in {\it Ecole d'été en Géométrie Tropicale} in 2006 \cite{P-06}:
\begin{Pro}Let $V\subset (\mathbb{C}^*)^n$ be an algebraic hypersurface and $\gamma$ be the logarithmic Gauss map. If $Crit(\Log_{\mid V})$ (resp. $Crit(\Arg_{\mid V})$) denotes the set of critical points of the logarithmic (resp. argument) mapping,
then we have the following equalities:  $Crit(\Log_{\mid V}) = Crit(\Arg_{\mid V}) = \gamma^{-1}(\mathbb{RP}^{n-1})$.
\end{Pro}

\begin{prooof}Let $f$ be the polynomial defining the hypersurface $V$, $reg (V)$ be the set of regular points of $V$, and $\gamma$ the logarithmic Gauss map
\[
\begin{array}{ccccl}
\gamma&:&reg (V)&\longrightarrow&\mathbb{CP}^{n-1}\\
&&(z_1,\ldots ,z_n)&\longmapsto&(z_1\frac{\partial f}{\partial z_1}:\ldots :z_n\frac{\partial f}{\partial z_n})
\end{array}
\]
In local coordinates $(w_1,\ldots ,w_n) = (\Log z_1,\ldots ,\Log z_n)$ one has $z_j \frac{\partial}{\partial z_j} = \frac{\partial}{\partial w_j}$, so $\gamma (e^w)$ is the normal direction to the set $\{ W \mid f(e^W)=0 \}$. Put $w_j = x_j + i \theta_j$, and let $C = A+iB =(A_1+iB_1,\ldots ,A_n+iB_n)$ be  a complex normal at $W$. We know that the complex tangent plane is given by the set of vectors $C$ such that $<C,W> = 0$, which is equivalently to say that we have the following equations
\[
\left\{ \begin{array}{lll}
<A,X>&=&<B,\Theta >\\
<B,X>&=&-<A,\Theta >
\end{array}
\right.
\]
where $X= (x_1,\ldots ,x_n)$ and $\Theta = (\theta_1,\ldots ,\theta_n)$. The maps $W\mapsto X$ and $W\mapsto \Theta$ are submersion if $A,\, B$ are linearly independent over $\mathbb{R}$. On the other hand, if $B=kA$ then we can see that we have the following
$$
<A, X> = k <A,\Theta> = -k <B, X> = -k^2 <A, X>
$$
which is equivalently to $<A,X> = <B,\Theta> = 0$.
\end{prooof}

\vspace{0.2cm}

\subsection{Complex tropical curves.} ${}$

\vspace{0.2cm}

Let $h$ be  a strictly positive real number and $H_h$ be the self
 diffeomorphism of $(\mathbb{C}^*)^2$ defined by :
\[
\begin{array}{ccccl}
H_h&:&(\mathbb{C}^*)^2&\longrightarrow&(\mathbb{C}^*)^2\\
&&(z_1,\ldots ,z_n)&\longmapsto&(\mid z\mid^h\frac{z}{\mid
  z\mid},\mid w\mid^h\frac{w}{\mid w\mid} ).
\end{array}
\]
which defines a new complex structure on $(\mathbb{C}^*)^2$
denoted by $J_h = (dH_h)\circ J\circ (dH_h)^{-1}$ where $J$ is the
standard complex structure.

\noindent A $J_h$-holomorphic curve $V_h$ is a curve
holomorphic with respect to the $J_h$ complex structure on
$(\mathbb{C}^*)^2$. It is equivalent to say that $V_h = H_h(V)$ where
$V\subset (\mathbb{C}^*)^2$ is an holomorphic curve for the
standard complex structure $J$ on $(\mathbb{C}^*)^2$.

\noindent Recall that the Hausdorff distance between two closed subsets $A,
B$ of a metric space $(E, d)$ is defined by:
$$
d_{\mathcal{H}}(A,B) = \max \{  \sup_{a\in A}d(a,B),\sup_{b\in B}d(A,b)\}.
$$
In this Note we take $E =\mathbb{R}^2\times (S^1)^2$, with the distance
defined as the product of the
Euclidean metric on $\mathbb{R}^2$ and the flat metric on $(S^1)^2$. We have the following definition
given by G. Mikhalkin in \cite{M2-04}.

\begin{Def} A complex tropical curve $V_{\infty}\subset
  (\mathbb{C}^*)^2$ is the limit (with respect to the Hausdorff
  metric on compact sets in $(\mathbb{C}^*)^2$) of a  sequence of a
  $J_h$-holomorphic curves $V_h\subset (\mathbb{C}^*)^2$ when
  $h$ tends to zero.
\end{Def}

\vspace{0.2cm}

\subsection{External hyperplanes and extra-pieces.}${}$

\vspace{0.2cm}

In this subsection we consider a complex algebraic hypersurface $V$ defined by a polynomial $f$ with Newton polytope $\Delta$. We denote by $\mathscr{A}$ the amoeba of $V$ and by $\Gamma$ its spine. Let us denote by $\tau$ the subdivision of $\Delta$ dual to $\Gamma$, and assume that $f$ is as $(1)$ (c.f.  \cite{RST-05} and \cite{PR-04} for more details).

\begin{Def}Let $\mathbb{R}^n$ be the universal covering of the real torus $(S^1)^n$. Let $\alpha$ and $\beta$ in the support of $f$. A hypersurface $H_{\alpha\beta}\subset \mathbb{R}^n$ is called {\em codual} (or corresponding) to an edge $E_{\alpha\beta}$ if it is given  by the following equation:
$$
\arg (a_{\alpha}) - \arg (a_{\beta}) + <\alpha - \beta , x> = \pi .
$$
In addition if $E_{\alpha\beta}$ is an external edge of $\tau$ (i.e., $E_{\alpha\beta}$ is a proper edge of the Newton polytope $\Delta$), then  $H_{\alpha\beta}$ is called an {\em external hyperplane}.
\end{Def}

\noindent We denote by $\Critv (\Arg )$ the set of critical values of the argument map.

\begin{Def} An {\em extra-piece} is a connected component $\mathcal{C}$ of $co\mathscr{A}\setminus \Critv (\Arg )$ such that the boundary of its closure $\partial \, \overline{\mathcal{C}}$ is not contained in the union of hyperplanes. 
\end{Def}

This means that its boundary contains at least one component (smooth) in the set of critical values of the argument map.

\section{Area Counted  with Multiplicity  of Coamoebas and the Alga map}

For a polynomial in two variables $f\in \mathbb{C}[z,w]$, let $V$ be the curve of its zero set in $(\mathbb{C}^*)^2$. We use the Alga map introduced by G. Mikhalkin and A. Okounkov in \cite{MO-07} defined as follow. Let $(\mathbb{Z}_2)^2$ be the subgroup of real points of the torus $(S^1)^2$, and  take now the quotient group $T = (S^1)^2/(\mathbb{Z}_2)^2$. Let us denote by $\Alg$ the composition of the argument map with the projection $pr : (S^1)^2\rightarrow T$. Let $\mathcal{R}\subset co\mathscr{A}$ be a subset of the regular values of the argument map. We define the {\em multiplicity} of $\mathcal{R}$ as the cardinality of the inverse image by the argument map of any point in $\mathcal{R}$.

\begin{Lem} Let $V$ be a complex algebraic plane curve defined by a polynomial $f$ with Newton polygon $\Delta$. Then the area of the coamoeba $co\mathscr{A}$ of $V$ counted with multiplicity and denoted by
$\Area_{mult} (co\mathscr{A})$ cannot exceed $2\pi^2 \Area (\Delta )$.
\end{Lem}
\begin{prooof} Let $\displaystyle{f(z,w) = \sum_{(i,j)\in \Delta\cap \mathbb{Z}^2} a_{ij}z^iw^j}$ be the polynomial defining the curve $V$.
 If the subdivision $\tau$ dual to the spine of the amoeba of the curve is a triangulation, then by  the deformation of the curve given by 
$\displaystyle{f_t(z,w) = \sum_{(i,j)\in \Delta\cap \mathbb{Z}^2} a_{ij}(et)^{\nu (i,j)}z^iw^j}$ with $t\in ]0; \frac{1}{e}]$ with $\nu (i,j) = - c_{(i,j)}$ the constant defined by Passare and Rullg\aa rd in \cite{PR-04} we have the result. Assume on the contrary  it is not the case. It is well known  that there exists  a polynomial $g$ with the same support as $f$ such that the arguments of its coefficients are the same as those of $f$ (see \cite{R-00}), and  the subdivision of $\Delta$ dual to the spine of the amoeba of the curve defined by $g$ is a triangulation. The result is immediate after deforming the coefficients of $g$ on those of $f$ (the deformation can be given by $(1-u)a_{\alpha ,\, g} + ua_{\alpha}$ with $u\in [0;1]$ and $a_{\alpha ,\, g}$ the coefficient of $g$ with index $\alpha$).

\end{prooof}

\vspace{0.2cm}

\begin{Lem} Let $V\subset (\mathbb{C}^*)^2$ be a no singular algebraic  curve defined by a polynomial $f$ such that $\mathbb{R}V$ is a Harnack curve. Then the coamoeba $co\mathscr{A}$ of $V$ is equal to the coamoeba of the complex tropical curve $V_{\infty ,\, f}$.
\end{Lem}

\begin{prooof} Let $f$ be the polynomial defining the curve $V$ which we assume  with real coefficients. Consider
$\displaystyle{\{ f_t\}_{t\in ]0; \frac{1}{e}]}}$ the family of polynomials $\displaystyle{\{ f_t\}}$
defined as above. We claim that the coamoebas of the curves $V_{f_t}$ for sufficiently small $t$ and those with $t$ sufficiently close to $\frac{1}{e}$ are the same. Indeed, let $V_{\infty ,\, f}$ be the complex tropical curve which is the limit when $t$ tends to zero of the $J_t$-holomorphic curves $H_t(V_{f_t})$ (abusive notation, because here we have precisely $H_h$  with $h=-\frac{1}{\log t}$ if we take the notation of 2.4), and $\Gamma_{\infty}$ be the tropical curve which is the limit of the amoebas $\mathscr{A}_{H_t(V_{f_t})}$. We know that the set of critical values of the argument map in our case is contained in the subgroup $(\mathbb{Z}_2)^2$ of the real points of the torus (recall that for a Harnack curves  we have 
$Crit_{\Arg} = \gamma^{-1}(\mathbb{RP}^1) = \mathbb{R}V$ see \cite{M1-00}). We can remark that the arguments of the coefficients of $f_t$ are equal to those of $f$, and then the coamoeba of the truncation of $f_t$ to any element of the subdivision $\tau$ is the same as that of the truncation of $f$ to the same element. Using Lemma 1.3 and 1.4 \cite{N4-08}, one can deduce that there was only one possibility.

\end{prooof}

\begin{Rem} ${}$
 \begin{itemize}
\item[(a)] \, Let $V\subset (\mathbb{C}^*)^2$ be a no singular algebraic  curve defined by a polynomial $f$ such that $\mathbb{R}V$ is a Harnack curve.  Let  $V_{\infty ,\, f}$ and $\Gamma_{\infty}$ as above and $v$ be a vertex of $\Gamma_{\infty}$, then  by the tropical localization we can see that $\Arg (\Log^{-1}(v)\cap V_{\infty ,\, f})$ is just two triangles of area $\pi^2/2$, and whose vertices are real and pairwise identified (they are real because those points depends on the coefficients of $f$, see section 3.6 in \cite{N1-06}, which are real, and there area is equal to $\pi^2/2$ because the linear part of the affine linear surjection sending the standard simplex to the simplex $\Delta_v$ dual to the vertex $v$ is in $SL_2(\mathbb{Z})$). Let $v_1$ and $v_2$ be two adjacent vertices of $\Gamma_{\infty}$ with $E_{v_1v_2}$ the edge between them. Let $U$ be a small open neighborhood of  the edge  $E_{v_1v_2}$, and $H_{v_1v_2}$ be the line in the torus corresponding to the edge $E^*_{v_1v_2}$ in the triangulation  of $\Delta$ dual to $E_{v_1v_2}$. 
By tropical localization (Viro's principle see \cite{V-90})  and  Kapranov's Theorem (see \cite{N1-06} and \cite{K2-00}), we claim that
the set $\overline{\Arg (V_{\infty ,\, f}\cap \Log^{-1} (U))}$ contains the line $H_{v_1v_2}$ in its interior. Indeed, on the contrary, this means that $H_{v_1v_2}$ is 
contained in the critical values of the argument map, which contradict the fact that the only critical values of the map  $\Arg$ are in the subset of real points of the torus $(S^1)^2$.
\begin{figure}[h!]
\begin{center}
\includegraphics[angle=0,width=0.4\textwidth]{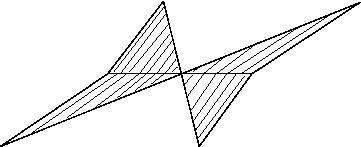}\quad
\includegraphics[angle=0,width=0.5\textwidth]{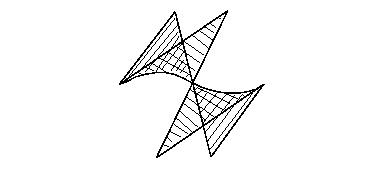}\quad
\includegraphics[angle=0,width=0.4\textwidth]{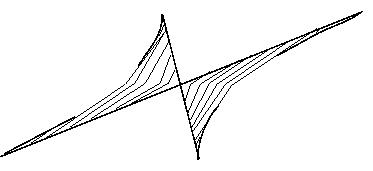}\quad
\caption{The first picture in the left is the Harnack case and the two others  are no-Harnack; extra-pieces and imaginary critical values of the argument map appear in the no-Harnack case.}
\label{}
\end{center}
\end{figure}

\item[(b)]\, Let ${\stackrel{\circ}{D}}_{std}$ be the interior of one triangle of the coamoeba of the standard line defined by $z+w+1=0$. If $\Delta '$ is a primitive simplex i.e., of area $\frac{1}{2}$, and $L$ is the linear part of the affine linear surjection of $\mathbb{R}^2$ sending the standard simplex to $\Delta '$, then for any $v\in (\mathbb{Z}_2)^2\subset (S^1)^2$ and $x,\, y$ two different points in ${\stackrel{\circ}{D}}_{std}$  we have $\tau_v\circ {}^tL^{-1}(x) \ne {}^tL^{-1}(y)$ where $\tau_v$ is the translation in the torus by the vector $v$. Indeed, assume that there exist two different points  $x$ and $y$ in ${\stackrel{\circ}{D}}_{std}$ and $v\in (\mathbb{Z}_2)^2$ such that $\tau_v\circ {}^tL^{-1}(x) = {}^tL^{-1}(y)$, then $x-y = {}^tL(v)$. So $x-y$ is a vector of the form $-\left( \begin{array}{l} k\pi\\ k'\pi \end{array}\right)$ which is impossible. Hence for any simplex $\Delta '$ of minimal area in $\mathbb{R}^2$, we have $\overline{pr\circ {}^tL_{\Delta '}^{-1}(co\mathscr{A}_{std})} = T$.

\end{itemize}
\end{Rem}

\section[Characterization of Coamoebas of Plane Curves Without Extra-Piece]{Characterization of Coamoebas of  Plane Curves Without Extra-Piece}

\begin{prooof}{\it of Theorem 1.1.} The implication (i) $\Rightarrow$ (ii) is obvious and the inverse is an immediate issue of Remarks 3.3 (a) and (b).
\end{prooof}

We say that a map is at most $p$-to-$1$ if the inverse image of any point in the target has at most $p$ points.

\begin{Lem} Let $V$ be a complex algebraic plane curve defined by a polynomial $f$  with Newton polygon $\Delta$, and let  $co\mathscr{A}$ be its coamoeba. Then if we have the equality $\Area_{mult}(co\mathscr{A}) = 2 \pi^2\Area (\Delta )$, then the map
$\Log_{\mid V} : V \rightarrow \mathbb{R}^2$ is at most $2$-to-$1$. 
\end{Lem}

\begin{prooof}  Let $\Gamma$ be the  spine of the amoeba of $V$ (viewed as tropical curve), and $\tau$ the subdivision of $\Delta$ dual to $\Gamma$. If we  take the same notation as in the proof of Lemma 3.2,  let $V_{\infty ,\, f}$ be the complex tropical curve equal to the limit of $H_t(V_{f_t})$ when $t$ tends to zero (with respect to the Hausdorff metric on compacts of $(\mathbb{C}^*)^2\approx \mathbb{R}^2\times (S^1)^2$ equipped with the product of the standard Euclidean metric and the flat metric on the torus).
Assume that the map $\Log_{\mid V} : V \rightarrow \mathbb{R}^2$ is not at most $2$-to-$1$. Let us denote by $F$ (resp. $F_c$) the set of critical values of the logarithmic (resp. argument) map, and we set $co\mathscr{A}\setminus F_c = \cup_{i=1}^l\mathcal{R}_i$. By Proposition 2.3, if $\Arg^{-1}(\mathcal{R}_i) = \cup_{j=1}^{n_i}\mathcal{R}_{ij}$ ($n_i$ is the multiplicity of the region $\mathcal{R}_i$) then for any $i$ and $j$ the set $\Log (\mathcal{R}_{ij})$ is a connected component of $\mathscr{A}\setminus F$. We have one of the two following cases: (i) the tropical curve $\Gamma$
has at least one vertex $v$ of multiplicity strictly greater than one (the multiplicity of a vertex $v$ is the area of a parallelogram spaned by two primitive vectors in the direction of any  two edges adjacent to $v$, see \cite{M2-04}); (ii) there exist at least two adjacent vertices $v_i$ and $v_j$  of $\Gamma$ such that if $U$ is a small open neighborhood of the edge $E_{v_1v_2}$ joining $v_1$ and $v_2$, then $\overline{\Arg (\Log^{-1}(U)\cap V_{})}$ don't contains a small neighborhood of the line $H_{v_1v_2}$ codual to the edge $E_{v_1v_2}$. In the  two cases (i) and (ii) we have $\Area_{mult}(co\mathscr{A}) < 2 \pi^2\Area (\Delta )$. Indeed, the first case means that the area of the coamoeba of the truncation of $f$ to the element $\Delta_v$ of $\tau$ dual to $v$ is such that $\Area_{mult}(co\mathscr{A}_{f^{\Delta_v}}) < 2 \pi^2\Area (\Delta_v)$, because in this case the truncation $f^{\Delta_v}$ is maximally sparse  or it has some coefficients which has no contribution on the spine of the amoeba  and then the coamoeba of the truncation has some extra-pieces. The second case means that we have extra-pieces and then $\Area_{mult}(co\mathscr{A}) < 2 \pi^2\Area (\Delta )$ (see Remark 3.3 (b)).
\end{prooof}

\begin{prooof}{\it of Theorem 1.2.} The implication (ii) $\Rightarrow$ (i) is obvious after Lemma 3.2 and Remark 3.3. (a), and the implication (i) $\Rightarrow$ (ii)  is an immediate issue of Lemma 4.3 and Theorem 1 in \cite{MR-00}.
\end{prooof}

\begin{figure}[h!]
\begin{center}
\includegraphics[angle=0,width=0.4\textwidth]{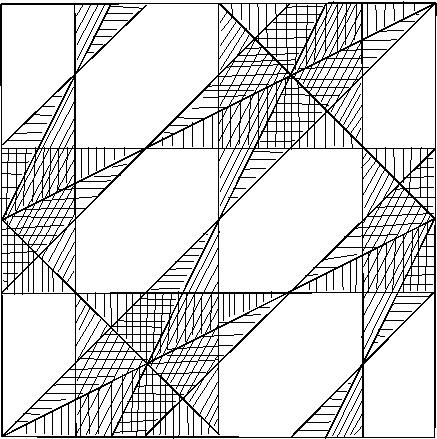}\quad
\caption{Coamoeba of a cubic with no extra-piece which is not Harnack and defined by a polynomial which is not maximally sparse.}
\label{}
\end{center}
\end{figure}

\end{document}